\documentclass[12pt]{article}
\usepackage{amssymb, amsmath, latexsym}
\usepackage{graphicx, setspace, amsthm, amsfonts, fancyhdr,fullpage}
\usepackage[all,2cell,ps]{xy}
	
\theoremstyle{plain}
\newtheorem{thm}{Theorem}[section]
\newtheorem{cor}[thm]{Corollary}

\newtheorem{lem}[thm]{Lemma}

\newtheorem{pro}[thm]{Proposition}

\theoremstyle{definition}
\newtheorem{dfn}[thm]{Definition}
\newtheorem{alg}[thm]{Algorithm}

\theoremstyle{remark}

\newtheorem{ntn}[thm]{Notation}
\begin{document}
\newcommand{\bm}[1]{\mbox{\boldmath{$#1$}}}
\newcommand{\ssection}[1]{
     \section[#1]{\sc #1}}  
\newcommand{\ssubsection}[1]{
     \subsection[#1]{\raggedright\it #1}}
\newcommand{\N}{\mathbb{Z}_{\geq 0}}
\newcommand{\Z}{\mathbb{Z}}

\numberwithin{equation}{subsection}

\title{Concrete and abstract structure of the
  sandpile group for thick trees with loops}
\author{William Chen and Travis Schedler}
\maketitle
\begin{abstract}
  We answer a question of L\'aszl\'o Babai concerning the abelian
  sandpile model.  Given a graph, the model yields a finite abelian
  group of recurrent configurations which is closely related to the
  combinatorial Laplacian of the graph. We explicitly describe the
  group elements and operations in the case of thick trees with
  loops---that is, graphs which are obtained from trees by setting
  arbitrary edge multiplicities and adding loops at vertices. We do
  this both concretely (by describing the so-called recurrent and
  identity configurations) and abstractly (by computing the group's
  abstract structure), and define maps identifying the two.
\end{abstract}

\ssection{Introduction} The \textit{abelian sandpile
  model}
(ASM) was defined in \cite{Dhar90} for any graph, generalizing 
the case of a grid from
\cite{BTW}.  This model is a prime example of \textit{self-organized
  criticality}
\cite{BTW} which has transformed the understanding of how complexity
arises in nature (see \cite{Turcotte} for a summary of applications to
modelling earthquakes, forest fires, landslides, etc.)  The ASM is 
related to the chip-firing game, introduced by Spencer
\cite{Spencer}, and later modified by Bj\"orner, Lov\'asz, and Shor
\cite{BLS}.  For general references on the abelian sandpile model, see
e.g.~\cite{BT, Dhar99}.

We will only consider undirected graphs.  An (undirected) abelian
sandpile model begins with any connected undirected graph (with
arbitrary edge multiplicities) with vertex set $V$ and edge set
$E$, together with a distinguished vertex $s \in V$ called the
\textit{sink}.  Then, a \textit{configuration} of the abelian sandpile
model is an element of $\N^{V_0}, V_0 := V \setminus \{s\}$, which
assigns a nonnegative integer to each vertex except the sink,
considered the number of grains of sand at that vertex.  A
configuration is \textit{stable} if the number of grains of
sand on each vertex is less than the degree of the vertex (otherwise
it is unstable).  Then, one obtains a dynamical system on the space of
configurations, where at each time step, each unstable vertex sends
one grain of sand along each incident edge.
Sending grains of sand to the
sink merely reduces the total number of grains.  Because every vertex has a
path to the sink, it turns out that this yields a well-defined
\textit{stabilization} of any configuration in finite time. The system
can remain in motion by continually adding grains (e.g., at random
vertices).

One basic question about any dynamical system is to find ``recurrent
configurations.''
For sandpiles, Dhar
\cite{Dhar90} defined these to be stable configurations that can be
obtained from any configuration by adding grains and stabilizing.
Dhar showed that they form an abelian group (under adding and
stabilizing), called ``the sandpile group,'' whose
abstract group structure is $\Z^n / \Delta \Z^n$, where $\Delta$ is
the matrix obtained from the Laplacian of the graph by deleting the
row and column corresponding to the sink.

However, the concrete structure (in terms of configuration on the
graph) of these recurrent configurations turns out to be quite
complicated and interesting.  For example, Creutz's paper \cite{Cre}
displays intricate fractal patterns of the identity element of this
group in the case of a square grid, which was further studied in
\cite{DM, LR}.

L\'aszl\'o Babai suggested that the authors study the sandpile group
for ``thick'' graphs: graphs which may have large edge multiplicities,
but whose combinatorial structure is simple.  In this note, we
consider \textit{thick trees with loops}: graphs which become trees
when all positive edge multiplicities are dropped to one, and all
loops are discarded.  We concretely describe the group elements and
operations in these cases, compute the abstract group structure, and
relate the two using recursive formulas.  To our knowledge, this is
the first result giving the complete, explicit structure of the
sandpile group in any class of undirected graphs.  In \cite{DRSV} (the
grid) and \cite{Bai, JNR, T} (other graphs), information about the
abstract group structure for more complicated graphs was given.

The note is organized as follows: in Section \ref{asm} we recall the
necessary background. In Section \ref{thtree} we state our main
theorem and its corollaries.  In Section \ref{pfs} we prove the theorem.

\ssection{The abelian sandpile model}\label{asm} We
review the theory of the (undirected) abelian sandpile model
(ASM), following \cite{BT}.

Let $A$ be a finite connected undirected graph with arbitrary edge
multiplicities, allowing loops.  Precisely, $A = (V,E,g)$, where $V$
is the set of vertices, $E$ is the set of edges, and
$g: E \rightarrow V^{(2)}$ is the incidence map.  Here $V^{(2)}$ is
the set of one- and two-element subsets of $V$ (i.e.,
unordered pairs of vertices).  We have $g(e)=\{i,j\}$ where $i,j$ are
the endpoints of $e$.  We require that $A$ be connected, i.e., there
is a path of edges connecting any two vertices.
\begin{dfn}
  An (undirected) \text{ambient space} is a collection $(A,s)$ where
  $A$ is an undirected, connected graph with vertex set $V$ and edge
  set $E$ (allowing multiple edges and loops), and $s \in V$ is a
  distinguished vertex, called the \textbf{sink}. All other vertices
  are called \textbf{ordinary}.  Let $V_0 := V \setminus \{s\}$ be the
  set of ordinary vertices.
\end{dfn}

\begin{ntn} We will use bold letters (e.g.~$\bm u$) to denote vectors,
  usually in $\N^{V_0}$.
\end{ntn}
\begin{dfn}
  A \textbf{configuration} of the ASM associated to $(A,s)$ is a vector
  $\bm u\in \N^{V_0}$, i.e., an assignment of a nonnegative integer
  $u_i$ to each ordinary vertex $i$ (called its \textbf{height}).
\end{dfn}

\begin{dfn}
  For any $i,j \in V$, define the \textbf{edge multiplicity}
  to be $e_{i,j} := \#\{a \in E: g(a) = \{i,j\}\}$. 
\end{dfn}

We may now define the notion of \textbf{stability} of a configuration:
\begin{dfn} Given a configuration $\bm u$, a vertex $i \in V_0$ is
  called \textbf{stable} if $u_i < \deg(i)$. Otherwise $i$ is called
  \textbf{unstable}.  A configuration is called \textbf{stable} if all
  vertices are stable.
\end{dfn}
\begin{dfn} \label{topdfn} If 
  $i \in V_0$ is unstable for $\bm u$, then define the state $\tau_i(\bm u)_i$
  by 
  $\tau_i(\bm u)_i := u_i - \deg(i) + e_{i,i}$ and
  $\tau_i(\bm u)_j := u_j + e_{i,j}$ for $j \neq i$.  Passing from
  $\bm u$ to $\tau_i(\bm u)$ is called ``toppling the vertex $i$.''
\end{dfn}
Given any unstable configuration, one may continually perform such
topplings $\tau_i$ until a stable configuration results, which we call
the \textbf{stabilization}.  
\begin{thm}\cite{Dhar90,BLS} (cf.~also \cite{Biggs99})
  Given any ambient space $(A,s)$ and any configuration $\bm u$ of its
  ASM, then there exists a unique stable configuration $\sigma(\bm
  u)$,
  which satisfies
  $\sigma(\bm u) = \tau_{i_1} \tau_{i_2} \cdots \tau_{i_m} (\bm u)$
  for some sequence of topplings
  $\tau_{i_1}, \tau_{i_2}, \ldots, \tau_{i_m}$.  
\end{thm}
(The theorem rests on an application of the Diamond Lemma,
since if two vertices are unstable at once, one may topple them
in either order with the same result.)
\begin{dfn}
  For each configuration $\bm u$, we denote by $\sigma(\bm u)$ its
  unique stabilization.
\end{dfn}
With this in mind, one may conclude that the set of stable configurations form
a monoid under vertexwise addition: 
\begin{dfn}\cite{Dhar90}
  Let $M$ be the set of stable configurations.
  Define the operation $\oplus$ on $M$ by
  $\bm u \oplus \bm v := \sigma(\bm u + \bm v)$, where $+$ is vector
  addition. 
\end{dfn}
\begin{lem}\cite{Dhar90}
  The operation $\oplus$ is commutative
  and associative, making $(M,\oplus)$ a commutative monoid, called
  the \textbf{sandpile monoid}.
\end{lem}

Even the abstract structure of the sandpile monoid (without discussing
concrete configurations) can be quite complicated. However, Dhar
\cite{Dhar90} noticed that the subset of ``recurrent'' configurations
has a much simpler abstract structure:
\begin{dfn}
  A configuration $\bm u$ is called \textbf{recurrent} if, for all
  configurations $\bm v$, there exists a configuration $\bm w$ such
  that $\bm v \oplus \bm w = \bm u$.
\end{dfn}
\begin{pro} \cite{Dhar90} The set of recurrent configurations, $G$, 
  is an abelian group under $\oplus$, called the \textbf{sandpile group}.  
\end{pro}
(Actually, as noticed in \cite{BT}, the above is true for any finite
commutative monoid, if one replaces $G$ by the unique minimal ideal.)
In \cite{Dhar90} (elaborated in \cite{DRSV}) there is
a surprisingly simple formula for the abstract structure of this group:
\begin{dfn}
  The \textbf{toppling matrix} 
  $\Delta=(\Delta_{ij})_{i,j \in V_0}$ is given by
\begin{equation}
\Delta_{ij} = \begin{cases} \deg(i) - e_{i,i}, & \text{ if }i=j, \\
                       -e_{i,j}, & \text{otherwise}.
\end{cases}
\end{equation}
\end{dfn}
Note that any two configurations which have the same stabilization are
related by an element in the row span of $\Delta$; also, adding a large
enough multiple of the sum of the rows $\Delta_i$ must give a unique
stabilization.  Using this, one has the
\begin{dfn} Let $\Lambda := \langle \bm \Delta_i \rangle_{i \in V_0}$ be
the lattice spanned by the rows $\bm \Delta_i$ of the toppling matrix.
\end{dfn}
\begin{thm}\cite{Dhar90,BT}
  The natural map $\N^{V_0} \twoheadrightarrow G$,
  $\bm c \mapsto \sigma(\bm c + \bm e^A)$ 
  descends to an isomorphism
  $\Z^{V_0}/ \Lambda \cong G$. In particular,
  $\#(G) = \det(\Delta)$.
\end{thm}
(Note: $G \cong \Z^{V_0}/ \Lambda$ and $\#(G) = \det(\Delta)$ are
due to \cite{Dhar90}; the map
$\bm c \mapsto \sigma(\bm c + \bm e^A)$ was pointed out in \cite{BT}.)

Dhar found an algorithm to test for recurrence, of linear
time in the combinatorial size of $A$ (the number of edges counted
without multiplicity plus the number of vertices):
\begin{dfn}
  Let $\bm \beta \in \N^{V_0}$ be given by
  $\bm \beta = \sum_{i \in V_0} \bm \Delta_i$. That is,
  $\bm \beta_i = e_{i,s}$.
\end{dfn}
We consider adding $\bm \beta$ to be ``toppling of the sink''. One has the
\begin{alg}[Burning algorithm]\label{ba2} \cite{Dhar90,MD92}
  To test whether a stable configuration ${\bm u}$ is recurrent,
  stabilize $\bm u + \bm \beta$.  Each vertex can topple at most once
  in the stabilization process, and $\bm u$ is recurrent iff every
  vertex topples (equivalently, iff $\bm u \oplus \bm \beta = \bm u$).
\end{alg}
The term ``burning algorithm'' (with which each vertex that
topples in the process is considered ``burned'') was introduced in
\cite{MD}.  The algorithm above was discovered by \cite{Dhar90}, and
the proof follows from \cite{Dhar90} and \cite{MD92}.  
\ssection{The
  sandpile group for thick trees}\label{thtree}
\begin{dfn}
  Following L. Babai, we call an undirected graph $T$ a
  \textbf{thick tree} if the underlying simple graph $T'$ obtained by
  reducing all nonzero edge multiplicities to 1 is a tree (i.e., a
  connected graph without cycles). A vertex $i$ is called a
  \textbf{leaf} if $\deg_{T'}i=1$ and $i\neq s$.
\end{dfn}
\begin{dfn}
  An undirected graph $T$ is called a \textbf{thick tree with loops}
  if, after removing all loops (dropping $e_{i,i}$ to zero for
  all $i \in V$), the obtained graph is a thick tree.
\end{dfn}
Recall that a graph is a tree iff, for any two vertices, there exists
a unique path from one to the other. Let $i, j \in V$.  We define a
partial ordering $\succeq$ on $V$, where $i \succeq j$ if $j$ is on
the unique path of $T'$ from $i$ to the sink. By definition, this
means that $j \succeq s$ for all $j \in V$.  Also, $i\succ j$ if the
condition $i \neq j$ is added. Vertex $p(j)$ is the \textbf{parent} of
$j$ if $p(j)\succ j$ and $p(j)$ is adjacent to $j$.

For a thick tree with loops, $T$, we define
$\succeq$ using the underlying tree $T'$. 

We have the following main result, whose proof is the content of
the next section.
\begin{ntn}\label{evntn}
  For any vertex $i \in V_0$, let $\bm \delta_i$ be the elementary
  vector given by $(\bm \delta_i)_j = \delta_{ij} = 1$ if $i = j$ and
  $0$ otherwise.
\end{ntn}
\begin{thm} \label{mt} The sandpile group
  $(G \subset \N^{V_0}, \oplus)$ for a thick tree $T$ with loops can
  be characterized as follows (recall Notation \ref{evntn}):
\begin{enumerate}
\item[(i)] $G = \{\bm u \in \N^{V_0}: \forall j \in V_0, 
  \deg(j) > u_j \ge \deg(j)-e_{j,p(j)} \}$.
\item[(ii)] The map
  $\phi: G \rightarrow \prod_{j \in V_0} \Z/e_{j,p(j)},$
  $\bm u \mapsto \sum_{j \in V_0} \bm \delta_j \bigl( \sum_{i \succeq
    j} u_i \pmod{e_{j,p(j)}} \bigr)$ is a group isomorphism. In particular,
  the abstract group structure does not depend on which vertex 
  is chosen to be the sink.
\item[(iii)] The map $\phi^{-1}$ may be expressed recursively as
  follows: For any $j \in V_0$,
  \begin{equation} \label{pif}
    \phi^{-1}(\bm v)_j = \begin{cases} \tilde v_j, & \text{if $j$ is a leaf}, \\
      \displaystyle f_j\Bigl(\tilde v_j - \sum_{i \succ j}
      \phi^{-1}(\bm v)_i\Bigr), & \text{otherwise,}
\end{cases}
\end{equation}
where for any $m \in \Z/p$, we let $\tilde m \in \{0,1,\ldots, p-1\}$
be its preferred representative, and for any $j \in V_0$,
$f_j: \Z \rightarrow \{\deg(j)-e_{j,p(j)}, \ldots, \deg(j)-1\}$ is
given (uniquely) by $f_j(m) \equiv m \pmod {e_{j,p(j)}}$.
\end{enumerate}
\end{thm}
\begin{cor} \label{asreccor} Let $\bm u$ be any configuration. Then
  the associated recurrent configuration $\sigma(\bm u + \bm e^T)$ is
  given recursively by
\begin{equation}
\sigma(\bm u + \bm e^T)_j=
\begin{cases} u_j - \lfloor u_j/e_{j,p(j)} \rfloor e_{j,p(j)} & 
      \mbox{ if }j\mbox{ is a leaf},\\ 
  u_j+\lambda_je_{j,p(j)}+\sum_{i\succ j}(u_i-\sigma({\bm u}+\bm e^T)_i)
      &\mbox{ otherwise},
\end{cases}
\end{equation}
where
$\lambda_j= \lfloor\bigl(\deg(j)-u_j-\sum_{i\succ j}(u_i-\sigma({\bm u
  + \bm e^T})_i)-1\bigr)/e_{j,p(j)}\rfloor$.
If $\sigma(\bm u)$ is already recurrent, we can delete the $\bm
e^T$'s.
\end{cor}
\begin{proof}
  Simply apply the map $\phi$ and then the map $\phi^{-1}$; we then
  expand $f_j$.
\end{proof}
\begin{cor}\label{thicktreeid}
  The identity configuration ${\bm e^T}$ 
is given by
\begin{equation}
  e^T_j=\begin{cases}0&\mbox{ if }j\mbox{ is a leaf},\\ \lambda_{j}e_{j, p(j)}- \sum_{i\succ j}e^T_{i}&\mbox{ otherwise},\end{cases}
\end{equation}
where
$\lambda_j=\lfloor \bigl(\deg(j)+\sum_{i\succ j}e^T_{i}-1\bigr)/e_{j,
  p(j)}\rfloor.$
\end{cor}
\begin{proof}
  Specialize Corollary \ref{asreccor} to the case of $\bm u = 0$.
  (Alternatively, specialize \eqref{pif} to the case of $\bm v = 0$
  and expand $f_j$.)
\end{proof}
\ssection{Proof of Theorem \ref{mt}} \label{pfs} We divide Theorem
\ref{mt} into Propositions \ref{mt1}, \ref{mt2}, and \ref{mt3},
corresponding to (refinements of) parts (i), (ii), and (iii),
respectively.
\begin{pro}\label{mt1}
In a thick tree $T$ with loops, the following statements are equivalent:
\item[(i)] Configuration ${\bm u}$ is recurrent.
\item[(ii)] For all $j\in V_0$, $\deg(j)>u_j \ge\deg(j)-e_{j,p(j)}$.
\end{pro}
\begin{proof}
  Take any stable configuration ${\bm u}$. We apply Algorithm
  \ref{ba2}.  This means that we try to stabilize $\bm u + \bm \beta$,
  and test whether or not every vertex topples.  Since $\bm \beta$
  adds sand only to the vertices adjacent to the sink, one easily sees
  that no vertex can topple before its parent.  We can view adding
  $\bm \beta$ as having ``the sink topple,'' and thus inductively, a
  configuration is recurrent iff every vertex can topple once its
  parent topples.  We conclude that every vertex fires iff, for every
  vertex $j$, the amount of chips needed to topple ($\deg(j) - u_j$)
  is less than or equal to the amount of chips it inductively receives
  from its parent ($e_{i,p(i)}$).
\end{proof}

\begin{lem}\label{lem: treemodequiv}
  The lattice $\Lambda$ can be
  characterized as
\begin{equation}
  \Lambda = \{\bm v:  e_{j,p(j)} \mid \sum_{i \succeq j} v_i, \forall j \in V_0\}.
\end{equation}
\end{lem}
\begin{proof}
  First, we show the inclusion $\subseteq$. It suffices to show that,
  if $\bm v = \bm \Delta_{k}$ for any $k \in V_0$, then
  $e_{j,p(j)} \mid \sum_{i \succeq j} v_i$ for any $j \in V_0$.

  For $k = p(j)$, one has
\begin{equation}
  \sum_{i \succeq j} v_i = -e_{j,p(j)} \equiv 0 \pmod{e_{j,p(j)}}.
\end{equation}
For $k = j$, one has
\begin{equation}
\sum_{i \succeq j} v_i = e_{j,p(j)} \equiv 0 \pmod{e_{j,p(j)}}.
\end{equation}
For all $k \in V_0 \setminus \{j,p(j)\}$, one has
\begin{equation}
  \sum_{i \succeq j} v_i = 0.
\end{equation}
The result is proved.

Next, we show the inclusion $\supseteq$.  For any $j \in V_0$, consider the sum
\begin{equation}
  \sum_{i \succeq j} \bm \Delta_i = -e_{j,p(j)} (\bm \delta_{p(j)} - \bm \delta_j).
\end{equation}
A simple inductive argument on $V_0$ shows that these elements span
the desired space
$W = \{\bm v: e_{j,p(j)} \mid \sum_{i \succeq j} v_i, \forall j
\in V_0\}$.
To do this, note that if $j$ is a leaf,
$\Delta_{jj} = e_{j,p(j)}$, and one may express all elements of
$W$ as a sum of such $\bm \Delta_j$'s and vectors which are zero in
the leaf components.
\end{proof}

As a consequence, we deduce the
\begin{pro} \label{mt2} One has an isomorphism
  $\bar \phi: \Z/\Lambda \cong \prod_{j
    \in V_0} \Z/e_{j,p(j)}$ which descends from the map
\begin{equation}
  \phi: \N^{V_0} \rightarrow \prod_{j \in V_0} \Z/e_{j,p(j)}, \quad
  \phi(\bm u)_{j} = \sum_{i \succeq j} u_i \pmod{e_{j,p(j)}}.
\end{equation}
Hence, the restriction
$\phi \bigl|_{G}: G \rightarrow \prod_{j \in V_0} \Z/e_{j,p(j)}$ is a
group isomorphism, and one has
$\phi(\bm c) = \phi \bigl|_{G} \; \circ \; \sigma(\bm c + \bm e^T)$,
where $\bm e^T \in G$ is the identity element, and
$\bm c \in \N^{V_0}$ is arbitrary.
\end{pro}

It remains to combine the two propositions to give an explicit inverse
of the isomorphism $\phi \bigl|_{G}$ (to express the sandpile group
$G \subseteq M \subseteq \N^{V_0}$ in terms of configurations), which
is part (iii) of Theorem \ref{mt}.
\begin{pro} \label{mt3} The inverse of $\phi \bigl|_G$ is given by
  \eqref{pif}.
\end{pro}
\begin{proof}
  Since $\phi \bigl|_G$ is an isomorphism and hence a bijection, we
  need only check that (i) the formula \eqref{pif} indeed gives a
  recurrent element (i.e., an element of $G$) for every
  $\bm v \in \prod_{j \in V_0} \Z/e_{j,p(j)}$, and (ii) the
  composition $\phi \circ \phi^{-1}$ is the identity.

  The first part (i) is obviously true by definition, using
  Proposition \ref{mt1} (we defined $f_j$ to satisfy this).  Let us
  check (ii).  To avoid confusion, let us denote the map recursively
  defined in \eqref{pif} by $(\phi^{-1})'$ for now (so replace every
  instance of $\phi^{-1}$ with $(\phi^{-1})'$. We then compute
\begin{equation}
  \phi \circ (\phi^{-1})'(\bm v)_j = f_j(v_j) - \sum_{i \succ j} 
  (\phi^{-1})'(\bm v)_i + \sum_{i \succ j} (\phi^{-1})'(\bm v)_i 
  \pmod{e_{j,p(j)}} = v_j.
\end{equation}
\end{proof}

\ssection{Acknowledgements} The authors would like to thank L\'aszl\'o
Babai for introducing them to the sandpile model through lectures at
the 2005 University of Chicago REU, supported by the University's
VIGRE grant; in particular for explaining all of the material in
Section 2, based on \cite{BT}.  They are also grateful to Babai for
proposing the problem, for providing many references, and for
extensive comments on earlier versions of this note.  The authors are
grateful to the participants and organizers of the aforementioned REU,
without whom this research would not have occurred.  The work of the
second author was partially supported by an NSF GRF.

\begin{singlespacing}

\noindent \ \\
\footnotesize{
{\bf W.C.}: California Institute of Technology, MSC 176, Pasadena, CA 91126, USA;\\
\hphantom{x}\quad\, {\tt chenw@caltech.edu}} \\
\footnotesize{
{\bf T.S.}: Department of Mathematics, University of Chicago, 
Chicago, IL
60637, USA;\\ 
\hphantom{x}\quad\, {\tt trasched@math.uchicago.edu}}
\end{singlespacing}
\end{document}